\documentclass{ifacconf}

\usepackage{amsfonts,amssymb, amsmath}
\usepackage{color}
\usepackage{graphicx}      
\usepackage{natbib}        
\usepackage[locale=DE]{siunitx}

\begin{document}
%
%
%
	
\begin{frontmatter}

\title{Feedback stabilization of forming processes\thanksref{footnoteinfo}} 

\thanks[footnoteinfo]{This work has been supported by DFG HE5386/19-1 and BA4253/11-1 of the Priority program 2183 on 'Property-Controlled Forming Processes'. }

\author[First]{Markus Bambach}  
\author[Second]{Michael Herty}
\author[Third]{Muhammad Imran}

\address[First]{Advanced Manufacturing, ETH Zurich, Switzerland 
	 (e-mail: mbambach@ethz.ch).}
\address[Second]{Institute of Applied Mathematics, RWTH Aachen University, Germany (e-mail: gerster@igpm.rwth-aachen.de)}
\address[Third]{Lehrstuhl Konstruktion und Fertigung, BTU Cottbus-Senftenberg, Germany (e-mail: muhammad.imran@b-tu.de)}

\begin{abstract}                
We are interested in the control of forming processes for nonlinear material models. To develop an online   control we derive a novel feedback law and prove  a  stabilization result. The derivation of the feedback control law is based on a Laypunov analysis of the time-dependent viscoplastic material models. The derivation uses the structure of the underlying partial differential equation for the design of the feedback control. Analytically, exponential decay of the time evolution of perturbations to desired stress--strain states is shown. We test the new control law numerically by coupling it to a finite element simulation of a deformation process.
\end{abstract}


\end{frontmatter}

\section{Introduction}

Metal forming processes allow for thermomechanical treatment and are thus able to convert the microstructure from an initial cast state into a desired state. Typically, this requires to use appropriate and a priori computed control inputs. For titanium alloys, \cite{BAMBACH2019249} propose punch velocity curves derived from the solution to an optimal control problem. However, even if the forming process is computed a priori uncertainties in the process and the material behavior may change the desired hidden microstructure and may lead to deviations from the desired material state. Model predictive control (MPC) based on a mean field microstructure model has been proposed by \cite{bambach2019}, but due to the non-linearity of the process, only a very much simplified cylinder compression test was considered, which can be treated using analytical process models. 
To overcome the limits of MPC, we aim at developing provable stabilizing feedback controls. For their derivation we consider a generic time--dependent viscoplastic nonlinear model and show that, at least in the spatially one-dimensional case, exponential damping of perturbations is theoretically possible. The obtained control relation is then applied to a finite-element simulation of a titanium forming process with strain hardening and softening in three spatial dimensions. 
\medskip
Mathematically, the derivation of the stabilizing feedback law follows recent results on stabilizing closed--loop control laws developed for 2$\times$2 hyperbolic
systems \cite{O1}. Examples of such dynamics describe the temporal and spatial evolution of gas \cite{G4} and water flow \cite{Prieur2008ab}. A key observation is that  also nonlinear viscoplastic materials are described by hyperbolic balance laws. Therein, suitable feedback boundary control laws have been derived. The  underlying tool for the study of the stabilization effect  are suitable Lyapunov functions measuring the deviation of the current state from the desired. Exponential decay of a continuous Lyapunov function under so-called \emph{dissipative} boundary conditions has been proven for example in ~\cite{L1,L4}.
For linearized balance laws explicit decay rates as well as numerical schemes have been presented e.g. in ~\cite{D1,Gerster2019,Knapp2019}. However, it is also known that if  the destabilizing effect of the source term is sufficiently large, the system cannot be controlled by this kind of feedback control~\cite{Bastin2011}.

\section{Derivation of Feedback Control Law}

As outlined in the previous section the derivation will be conducted in the spatially one-dimensional domain. For simplicity we assume that the domain is parameterized by $x \in [0,L]$ and we assume boundary conditions can be applied at $x=0$ and $x=L,$ respectively.  Since we are interested in stabilizing a desired state we denote by  $t \geq 0$ the time. Note that this is the time scale in which we want to stabilize the system, not the time scale under which  the forging process evolves. The time $t$ will only be required to derive the feedback law and to show stabilization of the dynamics.  On this scale we also assume that the material density $\rho$ is constant and normalized to $\rho = 1.$ Then, the evolution of the displacement $u=u(t,x)$, the displacement velocity $v(t,x)=\partial_t u(t,x),$  the stress $\sigma(t,x)$ and the total strain $\epsilon=\epsilon(t,x)$ are governed by the following balance and constitutive equations for $t \geq 0, x \in [0,1]$ 
\begin{align} 
	\partial_t v(t,x)  - \partial_x \sigma(t,x) &= 0, \\
	\partial_t \sigma(t,x)  &= \partial_t E \left( \epsilon - \epsilon^p \right)(t,x),
	\end{align}   
subject to initial conditions 
\begin{align}
	v(0,x) = v_0(x), \;  \sigma(0,x)=\sigma_0(x),
\end{align}
and boundary conditions to be specified below. In the previous equations $E>0$ is the elasticity modulus and $\epsilon^p$ the viscoplastic strain. The plastic part of the strain can be again a function of further variables, as the stress, the globularized volume fraction $X$ and the dislocation density $\bar{\rho}$. As seen below the precise form of $\epsilon^p$ is not required for the analysis. The detailed model for titanium alloy forming  is explained in the following section. 
\medskip Since $ \epsilon = \partial_x u(t,x)$ the previous system is rewritten as
\begin{align} 
	\partial_t v(t,x)  - \partial_x \sigma(t,x) &= 0, \\
	\partial_t \sigma(t,x)  - E \partial_x v(t,x) &= - E \partial_t \epsilon^p(t,x),
\end{align}   
The purpose is to stabilize within a time horizon $t\geq0$ the previous dynamics at a given state that yields the desired material properties. This state is assumed to be characterized by a stress $\sigma^*$ and a possibly non--constant displacement velocity $v^*(x).$ Perturbations to this state should be damped by a suitable feedback. The perturbation $U(t,x):=(\Delta v, \Delta \sigma)(t,x)$  where  $\Delta v(t,x):=v(t,x) - v^*(x)$ and $\Delta \sigma(t,x)=\sigma(t,x)-\sigma^*$ fulfills 
\begin{align} 
	\partial_t (\Delta v) - \partial_x (\Delta \sigma) &= 0, \\
	\partial_t (\Delta \sigma) - E \partial_x (\Delta v) =  - S^* \Delta \sigma. \label{Linearization}
\end{align}
Here, $S^* = \partial_\sigma \left( E \partial_t \epsilon^p(\sigma^*) \right)$ denotes the linearization of the plastic part of the strain at the desired state $\sigma^*.$ Depending on the dependence of the elastic part, the constant $S^*$ may contain additional terms due to the linearization with respect to globularized volume fraction or dislocation density.  
\medskip
The system is strictly hyperbolic and can be  diagonalized. The new variables $\mathcal{R}(t,x)=\left( \mathcal{R}^+, \mathcal{R}^- \right)(t,x) \in\mathbb{R}^2$ are given by  \begin{align}\label{trafo} 
	\mathcal{R} := T^{-1} U
\end{align} where $T$ contains eigenvectors of the matrix $A=\begin{pmatrix} 0 &  -1 \\ - E & 0 \end{pmatrix}.$ The eigenvalues of $A$ are  $\lambda_{1,2} = \pm \sqrt{E},$ respectively.  This shows that the system needs to be accompanied by a single boundary condition at $x=0$ and $x=L,$ respectively: 
\begin{align} \mathcal{R}^+(t,0) = K_0 \mathcal{R}^-(t,0), \; \mathcal{R}^-(t,L) = K_1 \mathcal{R}^+(t,L). 
 \end{align}
In the following we determine explicit characterization  on $K_0$ and $K_1$ such that $\| U(t,\cdot) \|^2_{L^2(0,L)}$ decays exponentially fast to zero over time. Regarding the previous conditions we note that those are local conditions, in the sense that $\mathcal{R}^+$ at $x=0$ is uniquely determined by $\mathcal{R}^-$ at $x=0$ and similarly at $x=L.$  We introduce the 
matrices 
\begin{align}   
	A = \begin{pmatrix} 0 & -1 \\ -E & 0  \\ \end{pmatrix}, \; 
S = \begin{pmatrix} 0 & 0 \\ 0 & S^*	\end{pmatrix},
\end{align}
and hence $T$ and $\Lambda = T^{-1} A T$ are given by 
\begin{align}   
	T:= \begin{pmatrix}
	-1 & 1 \\ \sqrt{E} & \sqrt{E}
\end{pmatrix}, \; 
	\Lambda  := \begin{pmatrix} \lambda_1 & 0 \\ 0 & \lambda_2	\end{pmatrix},
\end{align}
respectively. Further, we denote by $B:=T^{-1} S T.$ The previous transformations allow to rewrite the hyperbolic system as 
\begin{align} \label{LinearizedSystem}
	\partial_t \mathcal{R}(t,x) + \Lambda \partial_x \mathcal{R}(t,x) = -B \mathcal{R}(t,x), \\ 
	\begin{pmatrix} \mathcal{R}^+(t,0) \\ \mathcal{R}^-(t,L) \end{pmatrix}
	=
	\begin{pmatrix} & K_0 \\ K_1  \end{pmatrix}
	\begin{pmatrix} \mathcal{R}^+(t,L) \\ \mathcal{R}^-(t,0) \end{pmatrix}, \label{bcr}
\end{align}
and initial condition $\mathcal{R}(0,x)=T^{-1} U(0,x)$ for $x \in [0,L].$ 
For systems of the type \eqref{LinearizedSystem} general stabilization results exists, see e.g. \cite{Gerster2019}. However, we illustrate the main idea for its derivation based on the Lyapunov function, that is a weighted $L^2-$norm. This Lyapunov function has been introduced for hyperbolic systems in \cite{Coron2007aa}: 
\begin{align}\label{L}
		L(t) = \int_0^L  w^+(x;\hat{\mu}) \mathcal{R}^+(t,x)^2 + w^-(x;\hat{\mu}) \mathcal{R}^-(t,x)^2 dx,  
\end{align}
where for our purposes the weights $w^\pm$ are given by 
\begin{align}
		W(x;\hat{\mu}) := \begin{pmatrix} w^+(x;\hat{\mu}) & 0 \\ 0 & w^-(x;\hat{\mu})	\end{pmatrix}, \\ 
w^+(x;\hat{\mu}) = e^{-\frac{\hat{\mu}}{\sqrt{E}}x}, \; 
w^-(x;\hat{\mu}) = e^{-\frac{\hat{\mu}}{\sqrt{E}} (L-x)}.
 \end{align} 
We denote as $\lambda_{\min}(A)$ the minimal singular value of the matrix $A$ and as in \cite{Gerster2019} we define 
\begin{align}\label{w}
\mu (\hat{\mu}) := \hat{\mu} + \min\limits_{x\in[0,L]} \Big\{ \lambda_{\min} \Big[W(x;\hat{\mu}) B+B^T W(x;\hat{\mu})  \Big] \Big\}.
\end{align}
The main result is then as follows: Provided that there exists $\hat{\mu} \in \mathbb{R}^+$ such that 
 	\begin{align}
 	& \mu (\hat{\mu})
 	> 0, \label{ConditionsLyapunov1} \\
 	& \textup{exp}\Big( \hat{\mu} \frac{L}{2 \sqrt{E}} \Big) 
 	\max\big\{ |K_0|,|K_1|\big\} <1, \label{ConditionsLyapunov2}
 \end{align}
the Lyapunov function $L$ defined by equation \eqref{L} decays exponentially fast, i.e.,  
	\begin{align}\label{decay}
		L(t) \leq L(0) e^{-\mu (\hat{\mu}) t}.
		 \end{align}
Provided $\mathcal{R}$ is a differentiable solution to equation \eqref{LinearizedSystem} the equation \eqref{decay} is obtained by direct computation using \eqref{ConditionsLyapunov1} and \eqref{ConditionsLyapunov2}. 	 
	\medskip
Since $\mathcal{R}=T^{-1}U$ we obtain up to a constant 
\begin{align} \label{decay2}
C	\| U \|^2_{L^2(0,L)} \leq L(t) \leq \frac{1}C \| U \|^2_{L^2(0,L)}. 
	\end{align}	 
Hence, equation \eqref{decay} implies stabilization of the $L^2-$norm. Due to the linearity of the system \eqref{LinearizedSystem},  also higher-order norms can be stabilized as shown in \cite{Gerster2019}. Further note, that condition \eqref{ConditionsLyapunov2} is in fact a condition on the feedback constant $K_0$ and $K_1,$ respectively. This condition is called dissipativity condition in \cite{L4}.  
\medskip 
For the given system \eqref{LinearizedSystem} and any value $S^*$ the existence of $\hat{\mu}$ can be shown by the following estimate. In fact,  we evaluate \eqref{w} to obtain 
\begin{align*}
\mu (\hat{\mu})
	&=
	\hat{\mu}
	- \frac12
	\bigg| S^*
	\bigg| 
	\max\limits_{x\in[0,L]}
	\bigg\{
	3 \exp{-\frac{\hat{\mu} x}{\sqrt{E}}}
	+
	\exp{-\frac{\hat{\mu} (L-x)}{\sqrt{E}}}, \\ & 
	\exp{-\frac{\hat{\mu} x}{\sqrt{E}}}
	+
	3\exp{-\frac{\hat{\mu} (L-x)}{\sqrt{E}}}
	\bigg\} \\
	&\geq
	\hat{\mu} - 2\big| S^* \big|. 
\end{align*}
Hence, the decay rate $\mu(\hat{\mu} )$ is positive for example if 
\begin{align} \hat{\mu} := \big| S^* \big|. \end{align}
This shows that the precise modeling of the plastic strain as well as the desired state  part only enters in the estimate of the decay rate. In fact, an upper bound on $S^*$ would be sufficient to obtain $\hat{\mu}$ and implies condition \eqref{ConditionsLyapunov1}. The condition  \eqref{ConditionsLyapunov2} is fulfilled e.g. by the choice 
\begin{align}\label{formula}
K_0 = K_1 = 	\exp\left( 
		- \frac{L}{\sqrt{E}}
		| S^* |	\right). 
\end{align}
Boundary conditions in terms of $\mathcal{R}^\pm$ can be reformulated in terms of $U$ using equation \eqref{trafo}. By definition of $U$ we obtain by equation \eqref{bcr} 
\begin{align}
	\Delta v (t,0) &= - \frac{1 - K_0 }{ \sqrt{E} + K_0 \sqrt{E} } ( \Delta \sigma)(t,0), \\
	\Delta v (t,L) &=  \frac{ K_1 - 1 }{ \sqrt{E} +  K_1 \sqrt{E} } ( \Delta \sigma)(t,L).  
\end{align}
Applying the explicit form of the control \eqref{formula} yields the closed loop feedback law in terms of the desired state $(v^*(x), \sigma^*)$ as 
\begin{align}
v(t,0)
&=
v^*(0)
+
\frac{1}{\sqrt{E}}
\coth\bigg(
\frac{L}{\sqrt{E}}
\big| S^* \big|
\bigg)
\Big(
\sigma(t,0)-\sigma^*
\Big),
\end{align}
and similarly for $v(t,L).$ 


The purpose of the previous computations are to provide a framework for feedback control for general plastic material models. General stabilization problems  have been obtained and extended to $n\times n$ hyperbolic systems and we refer for example to \cite{O1} for more details. Also, we refer to \cite{Gerster2019} for a detailed analysis of conditions on the source term that are required for exponential stability as well as conditions on high--order norms.

\section{Numerical Results}

We apply the derived control law to a forming process of a titanium alloy. Here, a hybrid strain hardening and softening model is used that draws upon a physics-based dislocation density model to predict the flow stress $\sigma$  and globularized fraction $X$ as a function of strain rate $\dot{\epsilon}$ and temperature $T.$  The evolution 
of dislocation density $\bar{\rho}$ includes  generation of new dislocations as well as the annihilation of those \cite{GA1}. The evolution of $\bar{\rho}$ follows Taylor's equation and the the globularized fraction $X$ follows an Avrami-type model. This is used in the mixture rule to compute the flow stress when globularized and lamellar material volumes coexist. 
\medskip
The feedback law itself has been derived in terms of the displacement velocity $v$ at the boundary of the domain of the workpiece. 



\begin{figure}
 	\begin{center}
 		\includegraphics[width=8.4cm]{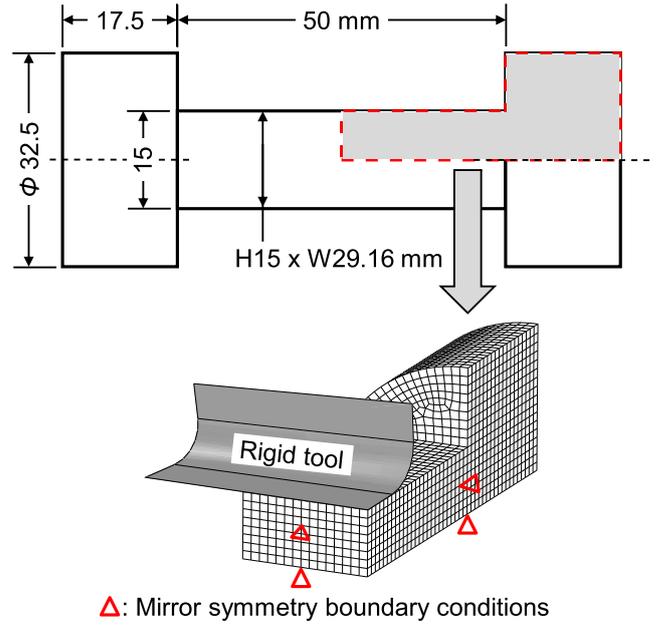}    
 		\caption{Workpiece geometry and the simplified finite element model for hot flat compression simulations.} 
 		\label{Figure1}
 	\end{center}
 \end{figure}

The derived feedback law is implemented in a finite element (FE) software.  The FE model of flat compression is set up in Abaqus Explicit software using the mirror-symmetry boundary conditions in the full scale model as shown in Fig. \ref{Figure1}. The reduced model is discretized using coupled temperature-displacement elements (C3D8T) with a mesh size of $1$ mm. The feedback law is applied at the boundary of the workpiece using the user amplitude subroutine (VUAMP) which defines the deformation velocity for the next time step based on the forced measured at the end of the previous time step. The feedback itself has been derived in terms of the velocity as a function of current stress $\sigma(t).$ However, the velocity boundary condition applied in the FE simulations has to be written in terms of the directly measured quantity, i.e. force $F=F(t,x)$ using the definition 
\begin{align}\label{F}
\sigma(t,x) = \frac{F(t,x)}A
\end{align}
 where $A$ is the cross-sectional area of the workpiece and $x$ is a point on the surface of the workpiece.

 The FE software provides the possibility of direct measurement of the actual force using “virtual sensor” in the history output requests. Then, the measurement of the virtual sensor can be directly called in the VUAMP. The material data for the simulations are taken from the hot compression tests of TNM-B1 at temperature $\ang{1150} C$ and $\ang{1200}  C$ and strain rates $0.0013, 0.005, 0.01$ and $0.05 \frac{1}s$ \cite{met9020220}. The specimen is deformed up to $3$ mm, so that the cross-section area remains constant. The simulations are performed using the initial velocity 
$$ v^*(0) = 1.5 \frac{mm}s,  A = 109.31 mm^2, L = 7.5 mm,$$
the  desired stress $\sigma^* = 146 MPa$ for $\ang{1150} C$
 and $68 MPa$ for $\ang{1200} C$, and material parameter $E = 9.2 \times 10^3 MPa$ for $\ang{1150} C$ and $E = 8.28 \times 10^3 MPa$ for $\ang{1200} C.$ 
 
\begin{figure}
	\begin{center}
		\includegraphics[width=8.4cm]{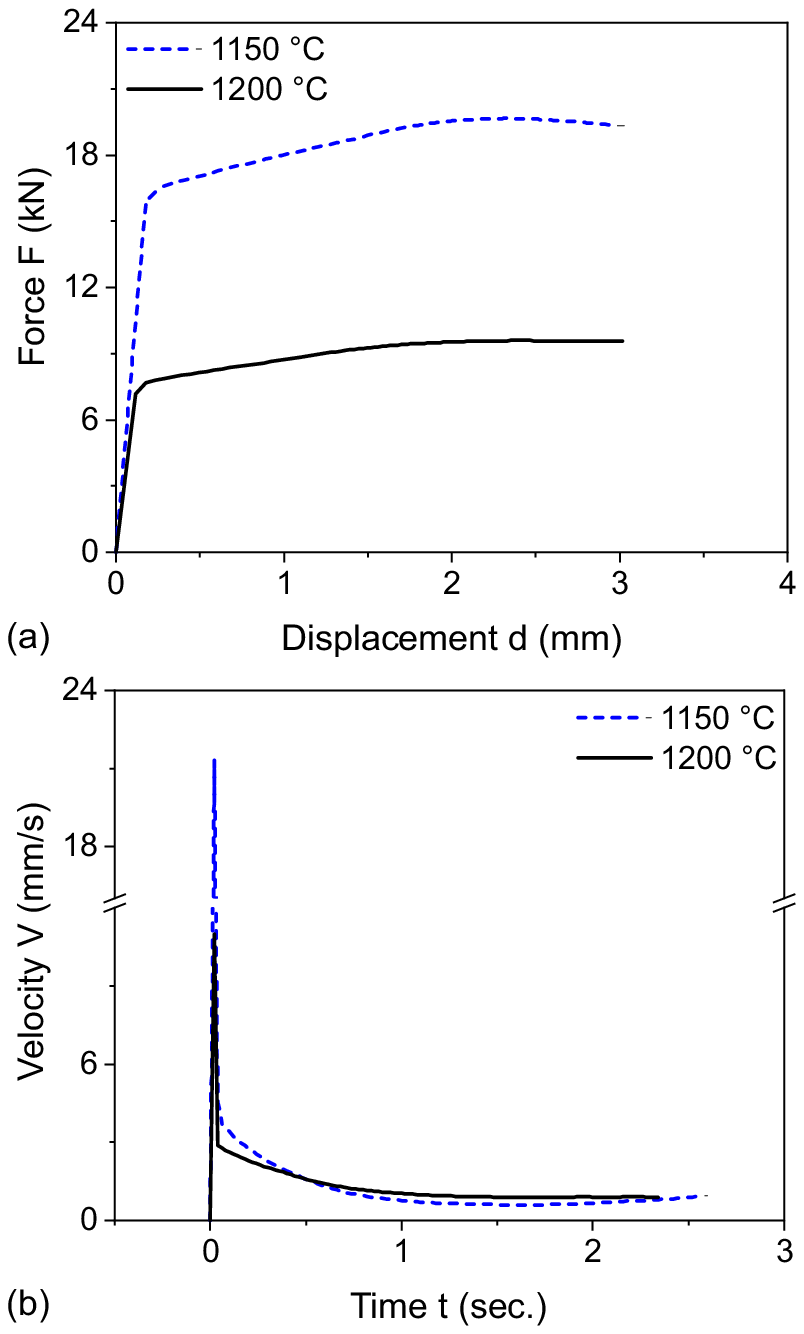}    
		\caption{Velocity-time profiles calculated using feedback control law (bottom) and resultant force-displacement curves (top) at temperatures $\ang{1150} C$ $(E = 9.2 \cdot 10^3 MPa)$ and $\ang{1200} C$ $(E = 8.28 \cdot 10^3 MPa).$} 
		\label{Figure2}
	\end{center}
\end{figure}

\begin{figure}

	\begin{center}
		\includegraphics[width=8.4cm]{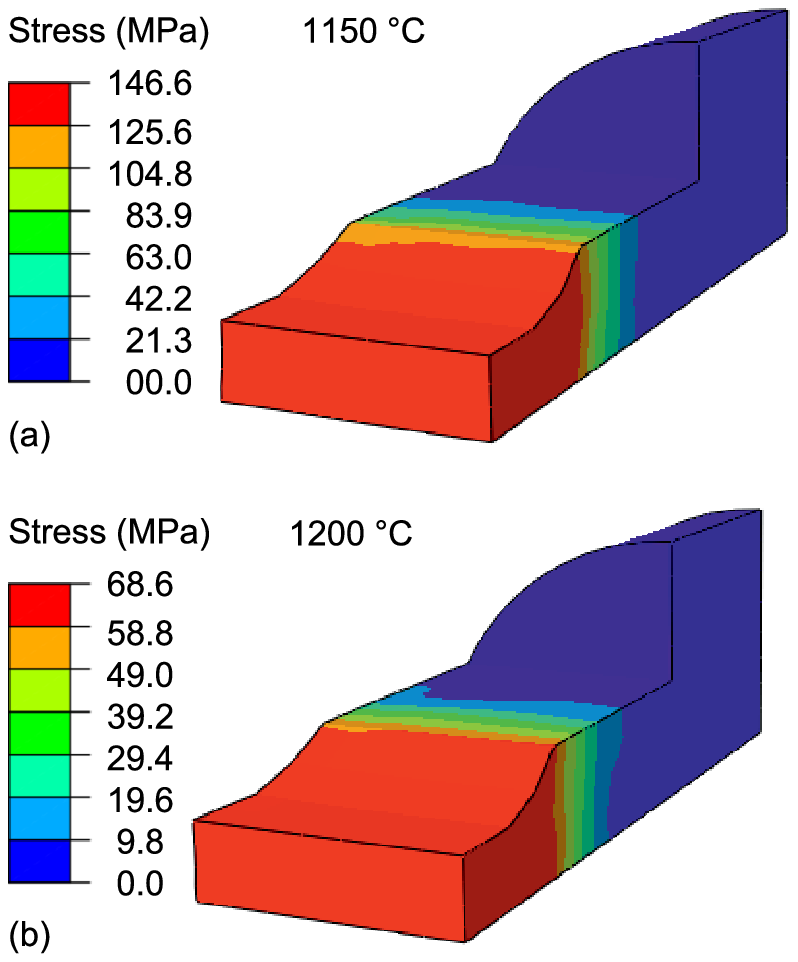}    
		\caption{Distribution of von-Mises stress in the deformed workpiece at temperature: (a) $\ang{1150} C$ $(E = 9.2 \cdot 10^3 MPa)$ and (b) $\ang{1200} C$ $( E = 8.28 \cdot 10^3 MPa)$.} 
		\label{Figure3}
	\end{center}
\end{figure}

\begin{figure}
	\begin{center}
		\includegraphics[width=5.4cm]{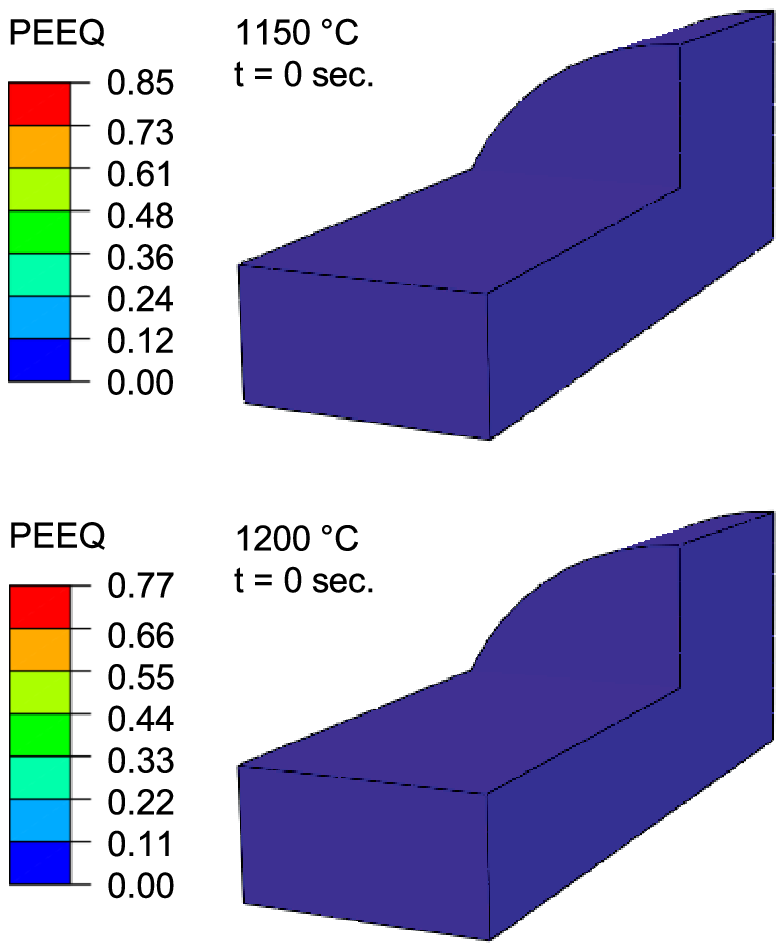} 
		\includegraphics[width=5.4cm]{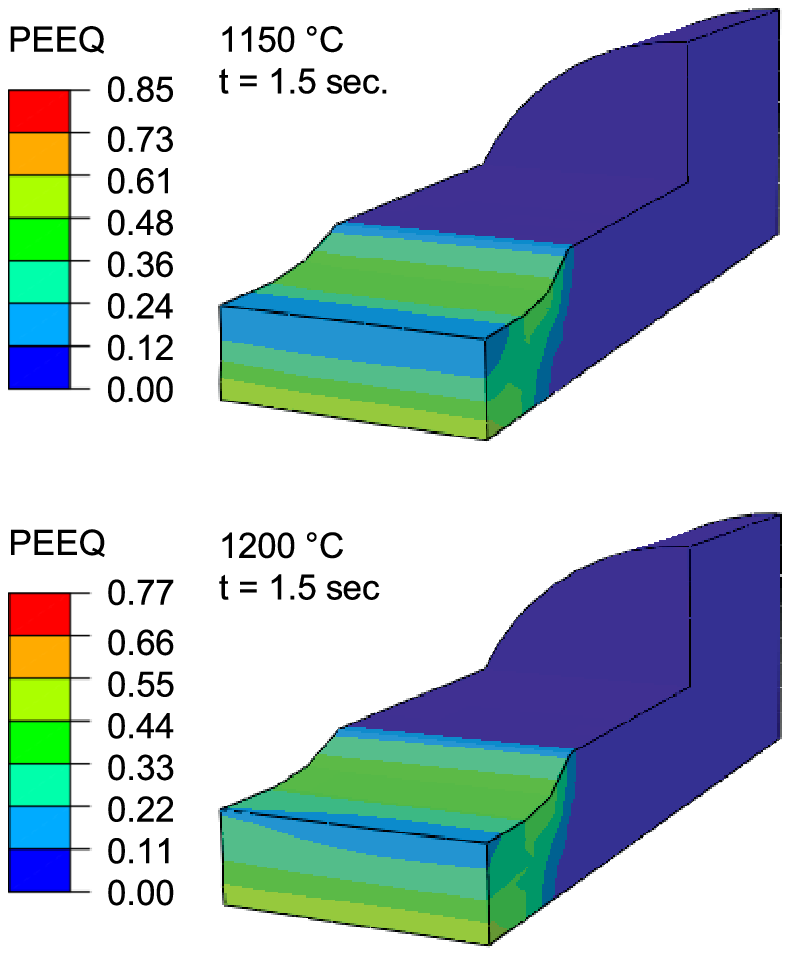} 
		\includegraphics[width=5.4cm]{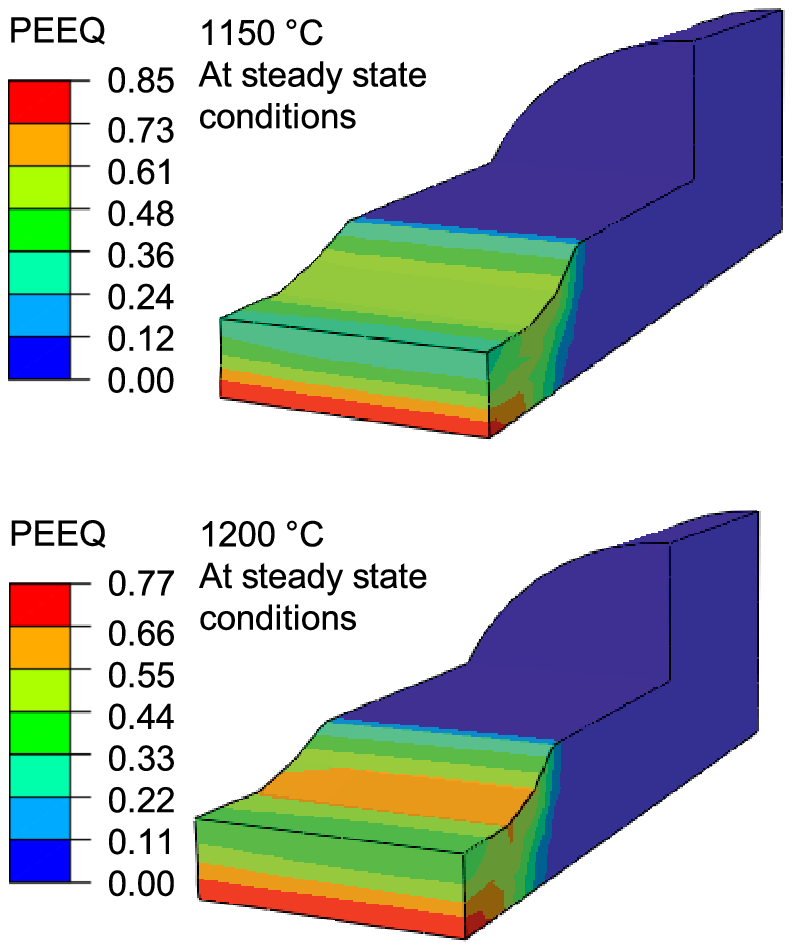}    
		\caption{Temporal evolution of the equivalent plastic strain (PEEQ) in the deformed workpiece at temperature:   $\ang{1150} C$ $(E = 9.2 \cdot 10^3 MPa)$ and $\ang{1200} C$ $( E = 8.28 \cdot 10^3 MPa)$ at time $0$ s, $1.5$ s and at steady state} 
		\label{Figure4}
	\end{center}
\end{figure}

The results of the feedback implementation in terms of force vs displacement and velocity decay over time are shown in Fig. \ref{Figure2}(a) and \ref{Figure2}(b), respectively. In the early stage of deformation, the difference between the force corresponding to desired stress and the actual force is large which leads to a rapid increase in the deformation velocity. With increasing deformation, the difference between the desired and the actual force decreases which causes a significant decrease in the deformation velocity. In the  steady-state condition of the force where the actual force meets the conditions of the desired force, the deformation velocity becomes equal to the initial velocity and the specimen deforms up to the required degree of deformation. At steady-state, the distribution of von--Mises stress and equivalent plastic strain (PEEQ) are presented in Fig. \ref{Figure3} and Fig. \ref{Figure4}, respectively.  For the plastic strain we also show the temporal evolution at two different temperatures to illustrate the effect of the control feedback. 

\medskip

The influence of the constants on model behavior is also analyzed. To achieve the desired stress in the workpiece with a fixed initial height, $v^*(0)$ and $E$ are two constants that can alter the model behavior. For the moment, $v^*(0)$ is kept constant. The constant $E$ is taken as $10 \times 10^3$ and $9\times 10^3$ MPa for $1150$ and $\ang{1200} C,$ respectively. The resultant material responses and the deformation velocities are compared in Fig. \ref{Figure6} (a) and (b), respectively. The results comparison shows that the material response is independent of constant $E.$ However, $E$ scales the deformation velocity and hence the process time. The independence of material response from E shows the robustness of the feedback control law where the process can be controlled with the desired stress and process time.

\begin{figure}
	\begin{center}
		\includegraphics[width=8.4cm]{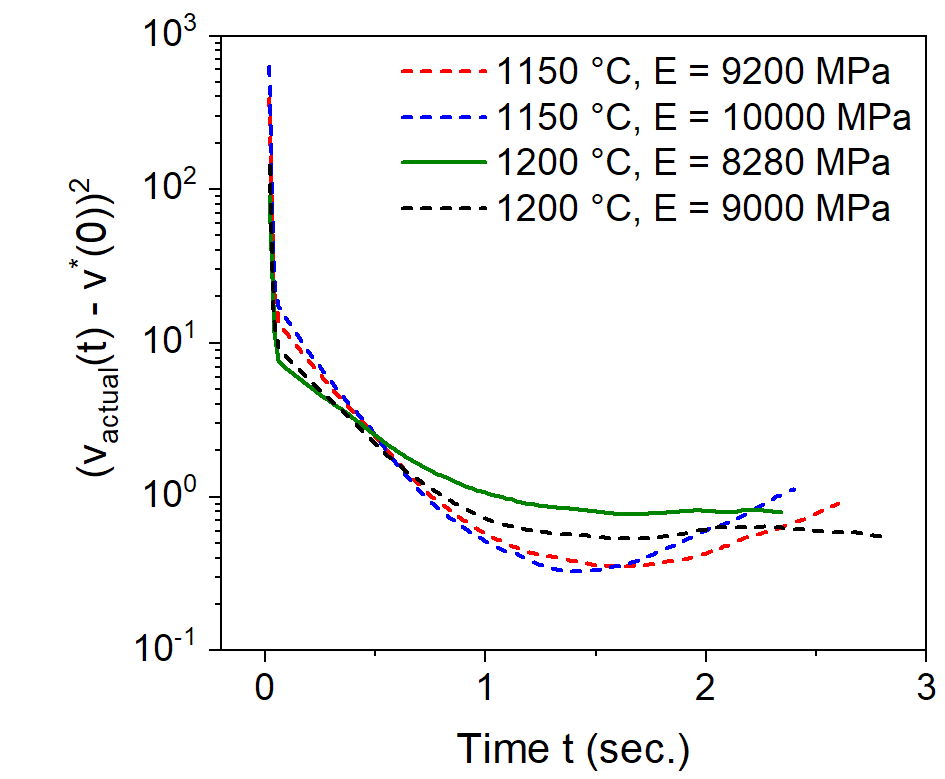}\\
		\includegraphics[width=8.4cm]{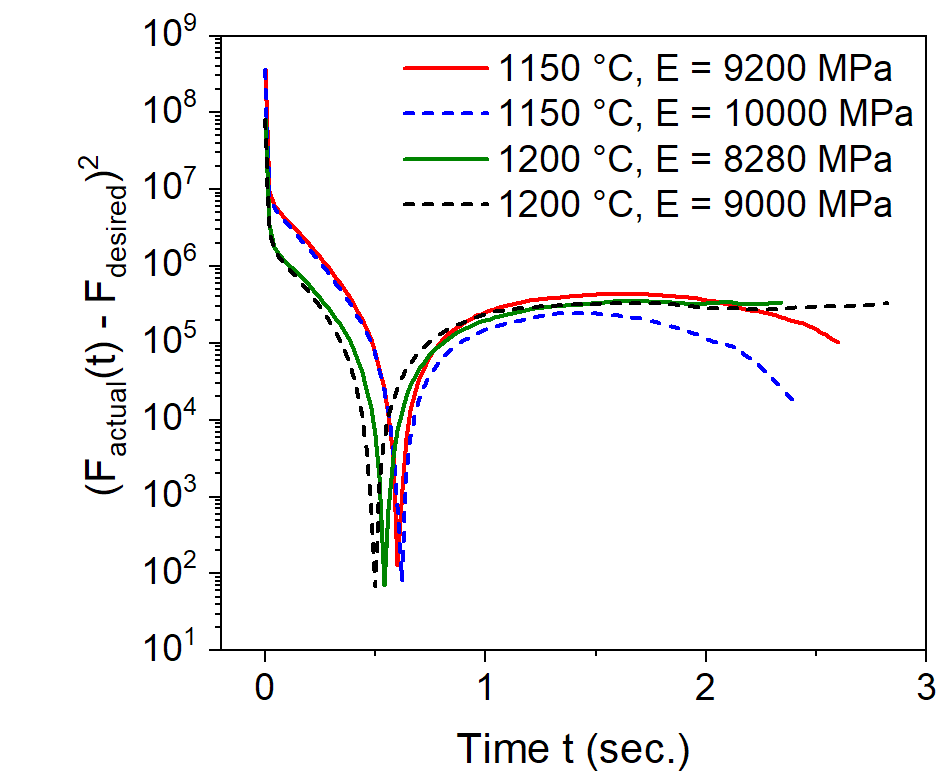} 
		\caption{Decay of the difference in stress $\sigma$ and $v$ at different temperatures in logarithmic scale.} 
		\label{Figure5}
	\end{center}
\end{figure}

\begin{figure}
	\begin{center}
		\includegraphics[width=8.4cm]{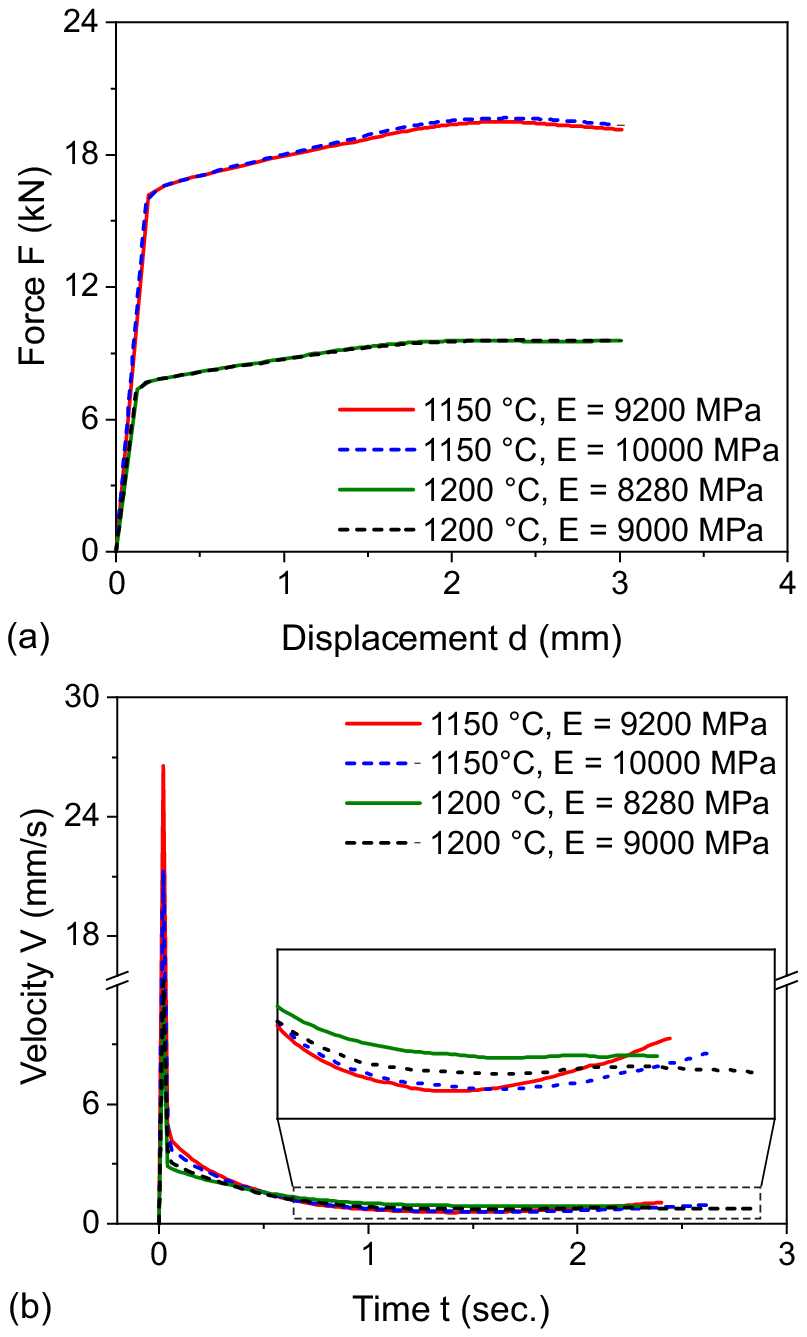}    
		\caption{	Velocity-time profiles calculated using feedback control law (bottom) and resultant force-displacement curves (top) at temperatures $\ang{1150} C$ $(E~=~10\cdot~10^3 MPa)$ and $\ang{1200} C$ $(E = 9 \cdot 10^3 MPa)$.} 
		\label{Figure6}
	\end{center}
\end{figure}

The theoretical estimate \eqref{decay2} suggests that we may expect exponential decay of $(v,\sigma)$ over time. Even so, the result is only in the spatially one-dimensional case we at least numerically investigate the decay in Fig. \ref{Figure5}. We observe decay in $v$ and in $F$ where $F(t)$ is computed along the boundary of the workpiece and includes the cumulative forces in the domain. Contrary to the one-dimensional result we observe a smaller decay rate in $F$. 

\section{Conclusion}

In the presented work we have shown a novel procedure to develop feedback control laws for nonlinear material models. Besides the derivation of the feedback law we show exponential stability for the proposed controlled system. Even so the derivation is limited to the spatially one-dimensional case we observe on realistic three-dimensional problems a similar performance. 

\begin{ack}
The authors thank the financial support of German Research Foundation (DFG)) through project BA4253/11-1 and HE5386/19-1 
of the Priority program 2183 on 'Property-Controlled Forming Processes'.
\end{ack}

\bibliography{main}             

\end{document}